\documentclass[11pt]{amsart}
\usepackage{amsmath}
\usepackage{amscd}
\usepackage{mathrsfs}
\usepackage{amssymb}

\newcommand{\R}{\mathbf{R}}
\newcommand{\Q}{\mathbf{Q}}
\newcommand{\Z}{\mathbf{Z}}
\newcommand{\C}{\mathbf{C}}

\newcommand{\set}[1]{\{ #1 \}}

\newcommand{\dbar}{\overline{\partial}}
\newcommand{\twopi}[1]{\frac{ #1 }{2\pi i}}
\newcommand{\Proj}{\mathbf{P}}
\newcommand{\Bl}[1]{\widetilde{#1}}
\newcommand{\ellip}[1]{\frac{\twopi{#1}\vartheta(\twopi{#1}-z)\vartheta'(0)}
{\vartheta(\twopi{#1})\vartheta(-z)}}

\newcommand{\coker}{\hbox{coker }}

\title[Rigidity and singular Chern numbers]{Rigidity of differential operators and Chern numbers of singular varieties}
\author{Robert Waelder}
\address{rwaelder@math.uic.edu}
\thanks{The author supported by NSF Post-doctoral Fellowship}

\begin{document}
\begin{abstract}
A differential operator $D$ commuting with an $S^1$-action is said to be rigid if the non-constant Fourier coefficients of $\ker D$ and $\coker D$ are the same. Somewhat surprisingly, the study of rigid differential operators turns out to be closely related to the problem of defining Chern numbers on singular varieties. This relationship comes into play when we make use of the rigidity properties of the complex elliptic genus--essentially an infinite-dimensional analogue of a Dirac operator. This paper is a survey of rigidity theorems related to the elliptic genus, and their applications to the construction ``singular" Chern numbers. 
\end{abstract}
\maketitle
\section{Rigidity of elliptic differential operators}
Let $D:\Gamma(E)\rightarrow \Gamma(F)$ be an elliptic operator maping sections of a vectorbundle $E$ to sections of $F$. If $D$ commutes with a $T=S^1$ action, then $\ker D$ and $\coker D$ are finite-dimensional $S^1$-modules. We define the character-valued index 
\begin{align*}
\hbox{Ind}_T(D) = \ker D- \coker D \in R(T)
\end{align*}
For example, if $D = d+d^* :\Omega^{\mathrm{even}}\rightarrow \Omega^{\mathrm{odd}}$ is the de Rham operator on a smooth manifold $X$ with a $T$ action, then by Hodge theory and homotopy invariance of de Rham cohomology, $\hbox{Ind}_T(D)$ is a trivial virtual $T$-module of rank equal to the Euler characteristic of $X$. In general, if $\mathrm{Ind}_T(D)$ is a trivial $T$-module, we say that $D$ is \emph{rigid}. In the case where $D$ is the de Rham operator, both $\ker D$ and $\coker D$ are independently trivial $T$-modules. However, more interesting cases exist where $D$ is rigid, but both $\ker D$ and $\coker D$ are nontrivial $T$-modules. For example, if $X$ is a spin manifold and $D:\Gamma(\Delta^+)\rightarrow \Gamma(\Delta^-)$ is the Dirac operator, then $D$ is rigid. It is instructive to sketch the proof of this fact, which is due to Atiyah and Hirzebruch \cite{AtiyahHirz}:

For simplicity, assume that $T$ acts on $X$ with isolated fixed points $\set{p}$, and that the action lifts to the spin bundles $\Delta^{\pm}$. At each fixed point $p$, $T_p X$ decomposes into a sum of one-dimensional complex representations of $T$ with weights $m_1(p),...,m_n(p)$, where $2n =\dim X$. If we view $\mathrm{Ind}_T(D)$ as a function of $t\in T$, then by the Lefshetz fixed point formula, 
\begin{align*}
\mathrm{Ind}_T(D) = \sum_p \frac{1}{\prod_{j=1}^n(t^{m_j/2}-t^{-m_j/2})}
\end{align*}
A priori, $\mathrm{Ind}_T(D)$ is a function only on the unit circle in $\C$. However, the above formula shows that we can analytically continue $\mathrm{Ind}_T(D)$ to a meromorphic function on $S^2$, with possible poles restricted to lie on the unit circle. But since $\mathrm{Ind}_T(D)$ is a virtual $T$-module, and therefore has a finite Fourier decomposition of the form $\mathrm{Ind}_T(D) = \sum a_n t^n$, all such poles on the unit circle must cancel. It follows that $\mathrm{Ind}_T(D)$ is constant. Furthermore, by taking the limit as $t\to \infty$, one sees that the character-valued index is identically zero. A similar proof shows that on a complex manifold, $\dbar+\dbar^*$ (whose corresponding index is the arithmetic genus) is rigid with respect to holomorphic torus actions.

The situation becomes more difficult if we investigate the rigidity of the twisted Dirac operators $D\otimes E$, where $E$ is an equivariant vectorbundle. For example, if $d_S=D\otimes (\Delta^+\oplus \Delta^-)$ is the signature operator on a spin manifold, the Lefschetz fixed point formula for the index of $d_S\otimes TX$ gives:
\begin{align*}
\mathrm{Ind}_T(d_S\otimes TX) = \sum_{p}\prod_{j=1}^n\frac{1+t^{-m_j(p)}}{1-t^{-m_j(p)}}\cdot\sum (t^{m_j(p)}+t^{-m_j(p)})
\end{align*}
Here $\pm m_j(p)$ are the weights of the $T$-action on the complexified tangent bundle of $X$ at $p$. The factors $\sum (t^{m_j(p)}+t^{-m_j(p)})$ come from the twisting of the rigid operator $d_S$ by $TX$. Thus, in this situation, the fixed point formula for $\mathrm{Ind}_T(d_S\otimes TX)$ has poles at $0$ and $\infty$, and we can no longer apply the same argument. 

It is therefore rather astonishing that, based on ideas from physics, Witten predicted the rigidity of an infinite sequence of twisted Dirac operators of this nature on a spin manifold. Witten's insight came from generalizing a quantum mechanics-inspired proof of the Atiyah-Hirzebruch theorem to its analogue in the setting of super string theory. We briefly sketch this point of view, as given in \cite{Witten2}: In super-symmetric quantum mechanics on a spin manifold $X$ (with one real fermion field), the Hilbert space of states corresponds to the space of square-integrable spinors. Quantization of the supercharge $Q_+$ yields the Dirac operator. In passing to super string theory, the Hilbert space of states should be interpreted as spinors on the loop space of $X$. It therefore makes sense to think of the quantization of the supercharge in this theory as a Dirac operator on the loop space. Now for any manifold $X$, the loop space of $X$ possesses a natural $S^1$ action given by rotating the loops. The fixed points of this action correspond to the space of constant loops, which we may identify with $X$ itself. Via formal application of the Atiyah-Bott-Lefschetz fixed point formula one can reduce the $S^1$ character-valued index of operators constructed out of $Q_+$ to integrals over $X$. To give an example, let $\Delta$ denote the spin bundle on the loop space. If we quantize a theory with two ferminionic fields $\psi_\pm$, the associated Hilbert space becomes $\Delta\otimes \Delta$. Now in finite dimensions, $\Delta\otimes\Delta$ corresponds to the de Rham complex. $\Delta\otimes\Delta$ therefore provides a good candidate for the de Rham complex on the loop space. At the classical level of this theory, one has an involution $\sigma$ on the space of superfields, sending $\psi_+\mapsto -\psi_+$ and $\psi_-\mapsto \psi_-$ which preserves the action-Lagrangian. When $X$ is spin, this involution descends to the quantum theory; the corresponding action of $\sigma$ on $\Delta\otimes \Delta$ may be interpreted as the Hodge star operator acting on forms. Consequently, one can construct out of $Q_+$ and $\sigma$ a canonical choice of a signature operator on the loop space. By the fixed point formula, its $S^1$-charactered valued index reduces to the index of
\begin{align*}
d_S\otimes\bigotimes_{n=1}^\infty \Lambda_{q^n}TX\otimes \bigotimes_{m=1}^\infty S_{q^m}TX = d_S\otimes\Theta_q
\end{align*}
over $X$. Here, for any vectorbundle $E$, we define $\Lambda_{q^m}(E) = 1+q^mE+q^{2m} E\wedge E+...$ and $S_{q^m}(E) = 1+q^m E+q^{2m}E^2+...$, where $q^m$ denote the weights of the induced $S^1$ action of an $S^1$-bundle over $X$. If $X$ itself has an $S^1$ action, the character-valued index of $d_S\otimes\Theta_q$ as a function of $e^{i\theta}$ may be interpreted as the signatures associated to a family of field theories parameterized by $\theta$. The rigidity of $d_S\otimes\Theta_q$ then follows from a formal application of deformation invariance of the index of Dirac operators on loop spaces. For details, see \cite{Witten1} or \cite{Witten2}. 

Note that since $d_S\otimes \Theta_q = d_S+2q d_S\otimes TX+...$, the rigidity of $d_S\otimes TX$ now follows from the rigidity of the $d_S\otimes \Theta_q$. It is interesting to point out that, although $d_S\otimes\Theta_q$ is defined on any oriented manifold, it is only rigid for spin manifolds. Heuristically this makes sense when we view $d_S\otimes \Theta_q$ as the signature operator on the loop space of $X$. For if $X$ is oriented, the signature operator $d_S$ is easily seen to be rigid. But the the loop space is oriented precisely when $X$ is spin.

Dirac operators on the loop space provide concrete examples of \emph{elliptic genera}. These are homomorphisms $\varphi:\Omega^{SO}\rightarrow R$ from the oriented cobordism ring to an auxiliary ring $R$, whose characteristic power series are defined in terms of certain elliptic integral expressions. 

The rigidity theorems of Witten were initially proven under restricted hypotheses by Landweber, Stong, and Ochanine \cite{LandweberStong,Ochanine}, and later proven in complete form by Bott, Taubes, and Liu \cite{BottTaubes,Liu}. The simplest and most direct proof was discovered by Liu, who observed that the modular properties of the elliptic genera implied their rigidity. We will provide a sketch of Liu's argument for the case of the complex elliptic genus, which is defined as the index of $\dbar\otimes E_{q,y}$ on an almost complex manifold, where $E_{q,y}$ is given by
\begin{align*}
E_{q,y}=y^{-n/2}\bigotimes_{n=1}^{\infty}\Lambda_{-yq^{n-1}}T''X\otimes\Lambda_{-yq^n}T'X\otimes\bigotimes_{m=1}^{\infty}S_{q^m}T''X\otimes S_{q^m}T'X
\end{align*}
Here $TX\otimes\C = T'X\oplus T''X$ denotes the decomposition of the complexified tangent bundle into holomorphic and anti-holomorphic components. By Riemann-Roch, the ordinary index of this operator is given by the integral
\begin{align*}
\int_X\prod_{T'X}\frac{x_j\vartheta(\twopi{x_j}-z,\tau)}{\vartheta(\twopi{x_j},\tau)}.
\end{align*}
Here $x_j$ denote the formal Chern roots of $T'X$, $y=e^{2\pi iz}$ and $q=e^{2\pi i\tau}$. $\vartheta(v,\tau)$ denotes the Jacobi theta function 
$$\vartheta(v,\tau) = \prod_{n=1}^\infty (1-q^n)\cdot q^{1/8} 2\sin \pi v \prod_{n=1}^\infty (1-q^n e^{2\pi iv})\prod_{n=1}^\infty (1-q^n e^{-2\pi iv})$$We will frequently refer to $\hbox{Ind}(\dbar\otimes E_{q,y})$ as $Ell(X;z,\tau)$.
The almost-complex version of Witten's rigidity theorem for this operator states that the complex elliptic genus of $X$ is rigid provided that $c_1(X) = 0$.

The idea of the proof is as follows: For simplicity, assume that the $T$-action on $X$ has isolated fixed points $\set{p}$, with equivariant weights $m_j(p)$ on $T'_p X$. By the Lefschetz fixed point formula,
\begin{align*}
\mathrm{Ind}_T(\dbar\otimes E_{q,y}) = \sum_p \prod_{j=1}^n\frac{\vartheta(m_j(p)u-z,\tau)}{\vartheta(m_j(p)u,\tau)}
\end{align*}

Write $\mathrm{Ind}_T(\dbar\otimes E_{q,y}) = F(u,z,\tau)$. 
It is evident from the fixed point formula that $F(u,z,\tau)$ is a meromorphic function on $\C\times\C\times\mathbf{H}$ which is holomorphic in $z$ and $\tau$. Let $z = \frac{1}{N}$ where $N$ is a common multiple of the weights $m_j(p)$. Then, using the translation formulas:
\begin{align*}
\vartheta(u+1,\tau) &= -\vartheta(u,\tau)\\
\vartheta(u+\tau,\tau) &= q^{-1/2}e^{-2\pi iu}\vartheta(u,\tau)
\end{align*}
it is easy to see that $F(u+1,\frac{1}{N},\tau) = F(u,\frac{1}{N},\tau)$ and that $F(u+N\tau,\frac{1}{N},\tau) = F(u,\frac{1}{N},\tau)$. Thus, for each fixed $\tau$, $F(u,\frac{1}{N},\tau)$ is a meromorphic function on the torus defined by the lattice $\Z\oplus N\Z\tau$. Suppose we could show that $F(u,\frac{1}{N},\tau)$ was in fact holomorphic in $u$. Then for each multiple $N$ of the weights $m_j(p)$ and for each $\tau\in \mathbf{H}$, $F(u,\frac{1}{N},\tau)$ would have to be constant in $u$. It would follow that $\frac{\partial}{\partial u}F(u,\frac{1}{N},\tau) \equiv 0$. Since this equation held for an infinite set of $(u,z,\tau)$ containing a limit point, it would hold for all $(u,z,\tau)$. Hence $F(u,z,\tau)$ would be independent of $u$, which is precisely the statement of rigidity for the operator $\mathrm{Ind}_T(\dbar\otimes E_{q,y})$. 

Thus, we are reduced to proving $F(u,z,\tau)$ is holomorphic. Let $\begin{pmatrix}a & b\\ c & d\end{pmatrix}\in SL_2(\Z)$ act on $\C\times\C\times\mathbf{H}$ by the rule $(u,z,\tau)\mapsto (\frac{u}{c\tau+d},\frac{z}{c\tau+d},\frac{a\tau+b}{c\tau+d})$. Using the transformation formula:
\begin{align*}
\vartheta(\frac{u}{c\tau+d},\frac{a\tau+b}{c\tau+d}) = \zeta(c\tau+d)^{\frac{1}{2}}e^{\frac{\pi i c u^2}{c\tau+d}}\vartheta(u,\tau)
\end{align*}
one sees that $F(\frac{u}{c\tau+d},\frac{z}{c\tau+d},\frac{a\tau+b}{c\tau+d})$ is equal to 
\begin{align*}
K\cdot \sum_p e^{-2\pi ic\sum_{j=1}^n m_j(p)u z/(c\tau+d)}\prod_{j=1}^n\frac{\vartheta(m_j(p)u-z,\tau)}{\vartheta(m_j(p)u,\tau)}
\end{align*}
where $K$ is a non-zero holomorphic function of $(u,z,\tau)$. Now the Calabi-Yau condition implies that the only possible $T$-action on $K_X$ is given by multiplication by a constant along the fibers. Since $\sum_{j=1}^n m_j(p)$ is the weight of the $T$-action induced on $K_X^*$, it follows that $\sum_{j=1}^n m_j(p)$ is the same constant for all $p$. We may therefore pull the expression $e^{-2\pi ic\sum_{j=1}^n m_j(p)u z/(c\tau+d)}$ outside of the above summation sign, and conclude that $F(\frac{u}{c\tau+d},\frac{z}{c\tau+d},\frac{a\tau+b}{c\tau+d}) = K' F(u,z,\tau)$, for $K'$ a non-zero holomorphic function.

Now the key observation: First, by the fixed point formula, $F(u,z,\tau)$ has possible poles only for $u = r+s\tau$, where $r,s\in \Q$. Moreover, since $F(u,z,\tau)$ is the character-valued index of an elliptic differential operator, the poles of $F(u,z,\tau)$ must cancel for $u\in \R$, since in that case $F(u,z,\tau)$ admits a Fourier decomposition $\sum b_m(z,\tau)e^{2\pi im u}$ (in a rigorous treatment of the subject, one must of course deal with convergence issues regarding this summation). Note that this is also the key observation in Bott and Taubes' proof. Thus, for $u$ a possible pole, write $u = \frac{n}{\ell}(c\tau+d)$, where $c$ and $d$ are relatively prime. By relative primality, we can find integers $a$ and $b$ so that $ad-bc = 1$, i.e., $\begin{pmatrix} a & b\\ c & d\end{pmatrix} \in SL_2(\Z)$. Then
\begin{align*}
K'\cdot F(\frac{n}{\ell}(c\tau+d),z,\tau) = F(\frac{n}{\ell},\frac{z}{c\tau+d},\frac{a\tau+b}{c\tau+d})
\end{align*}
where $K'\neq 0$. It follows that $F(u,z,\tau)$ is holorphic, which completes the proof.

The above rigidity theorem for the complex elliptic genus on a Calabi-Yau manifold has an interesting analogue for toric varieties, which has applications to the study of singular varieties. Let $\Sigma$ be a complete simplicial fan which corresponds to a smooth toric variety $X$. This means that $\Sigma$ is a finite union of cones $\set{C_i}$ inside the real vectorspace $N\otimes\R$, where $N$ is an integral lattice of rank $n$. For any two cones $C_1, C_2$ in $\Sigma$, we require that $C_1\cap C_2$ is a proper subcone, and that the union of the cones in $\Sigma$ covers all of $N\otimes\R$. The requirement that $\Sigma$ be simplicial simply means that the generators of $C_i$ are given by points in $N$. Recall that the data of $\Sigma$ gives rise to a natural scheme structure as follows: For each cone $C\subset\Sigma$, we define the sheaf of regular functions 
\begin{align*}
\Gamma(U_C) = \C[e^{f}]_{f\in S_C}
\end{align*}
where $S_C$ is the collection of linear functionals $f\in \mathrm{Hom}(N,\Z)$ that are positive along $C$. The underlying complex variety given by setting $U_C = \hbox{Specm }\Gamma(U_C)$ is, of course, the toric variety $X$. 

Note that inclusions of cones $C_1\subset C_2$ give rise to inclusions of open sets $U_{C_1}\subset U_{C_2}$. In particular, since every cone $C$ contains the point $0\in N$ as a subcone, every open set $U_C$ contains the open set 
\begin{align*}
U_0=\hbox{Specm }\C[e^{\mathrm{Hom}(N,\Z)}]\cong (\C^*)^n.
\end{align*}
The action of this complex torus on itself is easily seen to extend to all of $U_C$. In this way, $X$ inherits a natural action by a complex torus $T_\C$, with isolated fixed points.

There is a nice relationship between the $T_\C$-invariant divisors on a smooth toric variety and combinatorial data of its associated simplicial fan: the $T_\C$-invariant divisors on $X$ are in one-one correspondence with piecewise linear functionals on $\Sigma$. For example, if $f$ is a piecewise linear functional on $\Sigma$, then $f$ is completely determined by its values $f(v_i)$ on the generators $v_i\in N$ of the $1$-dimensional rays of $\Sigma$. These generators, in turn, define $T_\C$-Cartier divisors by the following prescription: If $C$ is a cone containing $v_i$, we define $\mathcal{O}(v_i)(U_C)=\Gamma(U_C)\cdot e^{v_i^*}$, where $v_i^*$ is the piecewise linear functional which is $1$ on $v_i$ and $0$ on the remaining $1$-dimensional rays of $\Sigma$. Otherwise, we set $\mathcal{O}(v_i)(U_C) = \Gamma(U_C)$. In this way, each piecewise linear $f$ gives rise to the divisor $D_f = \sum f(v_i)\mathcal{O}(v_i)$. In terms of this correspondence, it turns out there is a simple criterion for determining whether a $\Q$-divisor $D_f$ is linearly equivalent to zero: namely, $D_f\sim_{\Q} 0$ iff $f\in \hbox{Hom}(N,\Q)$.

Now, the canonical divisor $K_X = D_{f_{-1}}$, where $f_{-1}$ is the piecewise linear functional given by $f_{-1}(v_i) = -1$. Clearly if $\Sigma$ is complete, $f_{-1}$ cannot be given by a globally defined linear functional in $\hbox{Hom}(N,\Z)$. So compact smooth toric varieties are never Calabi-Yau, and consequently we can expect no rigidity properties for their complex elliptic genera. Note, however, that $TX$ is stably equivalent to $\bigoplus_{i=1}^\ell \mathcal{O}(v_i)$, where the sum is taken over all the $1$-dimensional rays $v_i$ of $\Sigma$. Thus, up to a normalization factor, the elliptic genus of $X$ is given by the index of $\dbar\otimes\xi$, where $\xi=$
\begin{align*}
\otimes_{i=1}^\ell\bigotimes_{n=1}^{\infty}\Lambda_{-yq^{n-1}}\mathcal{O}(v_i)^{-1}\otimes\Lambda_{-y^{-1}q^n}\mathcal{O}(v_i)\otimes
\bigotimes_{m=1}^{\infty}S_{q^m}\mathcal{O}(v_i)^{-1}\otimes S_{q^m}\mathcal{O}(v_i)
\end{align*}
We may view $\xi$ as a function of the $T_\C$-line bundle $\otimes_{i=1}^\ell \mathcal{O}(v_i)$. In this light, is natural to introduce, for any $T_\C$-line bundle $\otimes_{i=1}^\ell\mathcal{O}(v_i)^{a_i}$, with $a_i\neq 0$, the following vectorbundle, denoted as $\xi(a_1,...,a_\ell)$:
\begin{align*}
\otimes_{i=1}^\ell\bigotimes_{n=1}^{\infty}\Lambda_{-y^{a_i}q^{n-1}}\mathcal{O}(v_i)^{-1}\otimes\Lambda_{-y^{-a_i}q^n}\mathcal{O}(v_i)\otimes
\bigotimes_{m=1}^{\infty}S_{q^m}\mathcal{O}(v_i)^{-1}\otimes S_{q^m}\mathcal{O}(v_i)
\end{align*}
We may think of $\dbar\otimes\xi(a_1,...,a_\ell)$ as a kind of generalized elliptic genus for the toric variety $X$. The analogue of the Calabi-Yau condition for this generalized elliptic genus is the triviality of the $\Q$-line bundle $\otimes_{i=1}^\ell \mathcal{O}(v_i)^{a_i}$. In fact, if this bundle is trivial, then
\begin{align*}
\hbox{Ind}_T \dbar\otimes \xi(a_1,...,a_\ell) = 0 \in R(T)[[q,y,y^{-1}]]
\end{align*}
for any compact torus $T\subset T_{\C}$. To prove this, it suffices to assume that $T=S^1$ and that the $T$-action on $X$ has isolated fixed points. We can always find such a $T$ by first picking a dense $1$-parameter subgroup $\tau$ of a maximal compact subtorus of $T_\C$, and then letting $T$ be generated by a compact $1$-parameter subgroup whose initial tangent direction is sufficient close to that of $\tau$. Then the rigidity of $\dbar\otimes\xi(a_1,...,a_\ell)$ follows from Liu's modularity technique discussed above. To see that $\hbox{Ind}_T\dbar\otimes\xi(a_1,...,a_\ell)$ is identically $0$, we use the following trick observed by Hattori \cite{Hat}. Let $F(u,z,\tau) = \hbox{Ind}_T\dbar\otimes\xi(a_1,...,a_\ell)$. The modular properties of $F$ imply that $F(u+\tau,z,\tau) = e^{2\pi icz}F(u,z,\tau)$. Here $c$ is the weight of the $T$-action on the trivial bundle $\otimes_{i=1}^\ell \mathcal{O}(v_i)^{a_i}$. For a generic choice of $T\subset T_\C$, this weight will be non-zero. But since $F(u,z,\tau)$ is constant in $u$, we must have that $F(u,z,\tau) = e^{2\pi icz}F(u,z,\tau)$. This implies that $F(u,z,\tau) = 0$.

\section{Chern numbers of singular varieties}
We now turn to the problem of defining Chern numbers on singular varieties, a subject which at first glance appears unrelated to the discussion above. In what follows we will find that rigidity theorems provide a powerful tool in solving these types of problems. We first discuss some background.

If $X$ is a smooth compact almost-complex manifold of dimension $n$, the Chern numbers of $X$ are the numbers of the form
\begin{align*}
c_{i_1,...,i_n} = \int_X c_1^{i_1}\cdot c_2^{i_2}\cdots c_n^{i_n}
\end{align*} 
where $c_i$ denotes the $i$th Chern class of $T'X$ and $i_1+2i_2+...+ni_n = \dim X$ (so that the total degree of the integrand is $2n$). Chern numbers are easily seen to be functions on the complex cobordism ring $\Omega^*_U$. Moreover, they completely characterize $\Omega^*_U$ in the sense that two almost complex manifolds with the same Chern numbers must be complex cobordant.

Much of algebraic geometry consists of efforts to extend techniques from the theory of smooth manifolds to singular varieties. Minimal model theory suggests that one can approach this problem by working on a smooth (or ``nearly smooth'') birational model of a given singular variety $X$. For a special combination of Chern numbers, this approach works without any difficulties: namely, the Chern numbers defining the Todd genus. For if $X$ is a smooth complex manifold, the Todd genus of $X$ is given by the alternating sum $\chi_0(X) = \sum_{i=0}^n (-1)^i \dim H^{i,0}_{\dbar}(X)$. By Hartog's theorem, the space of holomorphic $i$-forms is birationally invariant, and is therefore well-defined even when $X$ is singular, by passing to a resolution of singularities. On the other hand, if $X$ is smooth, then by Riemann-Roch,
\begin{align*}
\chi_0(X) = \int_{X}\prod_{i=1}^{n}\frac{x_i}{1-e^{-x_i}}
\end{align*}
where $x_i$ denote the formal Chern roots of the holomorphic tangent bundle. The combination of Chern numbers obtained by performing the above integration therefore makes sense for any compact singular variety defined over $\C$.

More generally, we consider the following naive attempt at defining combinations of Chern numbers on $X$: Simply let $Y$ be a smooth minimal model of $X$ and define $c_{i_1,...,i_n}(X) = c_{i_1,...,i_n}(Y)$. Nevermind the implicit assumption that the minimal model program holds. The main problem with this approach is that we should not expect a unique choice of a minimal model $Y$. In general, $X$ will have various minimal models which differ from each other by codimension-$2$ surgeries called flips and flops. A priori, it is not at all clear what combinations of Chern numbers will be preserved under these operations.

In \cite{Totaro} Totaro set out to classify the combinations of Chern numbers invariant under classical flops. Here we say that two varieties $X_1$ and $X_2$ differ by a classical flop if they are the two small resolutions of an $n$-fold $Y$ whose singular locus is locally the product of a smooth $n-3$-fold $Z$ and the $3$-fold node $xy-zw = 0$. More precisely, $X_1$ and $X_2$ are constructed as follows: blowing up along $Z$ defines a resolution of $Y$ whose exceptional set is a $\Proj^1\times\Proj^1$ bundle over $Z$ with normal bundle $\mathcal{O}(-1,-1)$. Here $\mathcal{O}(-1,-1)$ denotes the line bundle over a $\Proj^1\times\Proj^1$-bundle which coincides with the tautological bundle along each $\Proj^1$ direction. Blowing down along either of these $\Proj^1$ fibers therefore produces two distinct small resolutions $X_1$ and $X_2$ of $Y$.

Totaro demonstrated that the combinations of Chern numbers invariant under classical flops were precisely the combinations of Chern numbers encoded by the complex elliptic genus in the Riemann-Roch formula. We sketch the first half of his argument--namely, that the complex elliptic genus is invariant under classical flops. As $X_1$ and $X_2$ are identical away from their exceptional sets, their difference $X_1-X_2$ is complex cobordant to a fibration $E$ over $Z$. In fact, if the exceptional sets of $X_1$ and $X_2$ are the $\Proj^1$-bundles $\Proj(A)$ and $\Proj(B)$ corresponding to the rank $2$ complex bundles $A$ and $B$ over $Z$, then as a differentiable manifold, $E$ is simply the $\Proj^3$ bundle $\Proj(A\oplus B^*)$ over $Z$. Now the way that $E$ is actually constructed is by taking a tubular neighborhood of $\Proj(A)\subset X_1$ and gluing it to a tubular neighborhood of $\Proj(B)\subset X_2$ along their common boundaries (which are both diffeomorphic to $Z\times S^3$). The crucial point is that the stably almost complex structure on $E$ induced by this construction makes $E$ into an $SU$-fibration. That is, $E$ is a $\Proj^3$-bundle whose the stable tangent bundle in the vertical direction has a complex structure satisfying $c_1 = 0$. He calls these fibers ``twisted projective space" $\Bl{\Proj}_{2,2}$. The fiber-integration formula implies that $Ell(E;z,\tau) = \int_Z Ell_T(\Bl{\Proj}_{2,2};z,\tau,x_1,...,x_4)\cdot\mathcal{E}ll(Z;z,\tau)$. Here $\mathcal{E}ll(Z;z,\tau)$ is the cohomology class which appears as the integrand in the Riemann-Roch formula for the elliptic genus of $Z$. More importantly, $Ell_T(\Bl{\Proj}_{2,2};z,\tau,x_1,...,x_4)$ denotes the character-valued elliptic genus of $\Bl{\Proj}_{2,2}$ with the standard $T^4$ action, with the generators $u_1,...,u_4$ of the Lie algebra of $T^4$ evaluated at the Chern roots $x_1,...,x_4$ of $A\oplus B$. Since $\Bl{\Proj}_{2,2}$ is an $SU$-manifold, by the Witten rigidity theorem, $Ell_T(\Bl{\Proj}_{2,2};z,\tau,x_1,...,x_4) = \hbox{const}$. Thus, the elliptic genus of $E$ is simply the product $Ell(\Bl{\Proj}_{2,2};z,\tau)\cdot Ell(Z;z,\tau)$. Moreoever, since $\Bl{\Proj}_{2,2}$ is cobordant to $Y_1-Y_2$, where $Y_i$ are the small resolutions of a $3$-fold node, and since classical flopping is symmetric for $3$-folds, $\Bl{\Proj}_{2,2}\sim Y_2-Y_1$. Hence $\Bl{\Proj}_{2,2}\sim 0$ in the complex cobordism ring. We therefore have that $Ell(X_1;z,\tau) - Ell(X_2;z,\tau)=Ell(\Bl{\Proj}_{2,2};z,\tau)\cdot Ell(Z;z,\tau)=0$.

An obvious consequence of the above discussion is that for varieties $Y$ whose singular locus is locally the product of a smooth variety with a $3$-fold node, it makes sense to define the elliptic genus of $Y$ to be the elliptic genus of one of its small resolutions. However, most singular varieties fail to possess even one small resolution. It is therefore natural to ask whether one can continue to define the elliptic genus for a more general class of singularities. The right approach to answering this question is to expand one's category to include pairs $(X,D)$, where $X$ is a variety and $D$ is a divisor on $X$ with the property that $K_X-D$ is $\Q$-Cartier. A map $f:(X,D)\rightarrow (Y,\Delta)$ in this category corresponds to a birational morphism $f:X\rightarrow Y$ satisfying $K_X-D = f^*(K_Y-\Delta)$. The idea is to first define the elliptic genus for smooth pairs $(X,D)$ in such a way that $Ell(X,D;z,\tau)$ becomes functorial with respect to morphisms of pairs. Given two resolutions $f_i:X_i\rightarrow Y$ of a singular variety $Y$, with $K_{X_i}-D_i = f^*K_Y$, we could then find resolutions $g_i:(M,D)\rightarrow (X,D_i)$ making the following diagram commute:
$$\begin{CD}
(M,D) @> g_1 >> (X_1,D_1) \\
@V g_2 VV       @VV f_1 V\\
(X_2,D_2) @> f_2 >> (Y,0)
\end{CD}$$ 
Functoriality of the elliptic genus would then imply that 
$$Ell(X_1,D_1;z,\tau)=Ell(M,D;z,\tau)=Ell(X_2,D_2;z,\tau).$$ It would then make sense to define 
$Ell(Y;z,\tau)\equiv Ell(X_1,D_1;z,\tau).$

 One can simplify this approach by making two observations. First, by introducing further blow-ups, one can always assume that the exceptional divisors $D_i\subset X_i$ have smooth components with simple normal crossings. (Such resolutions are called ``log resolutions".) Second, by a deep result of Wlodarczyk \cite{W}, the birational map $(X_1,D_1)\dashrightarrow (X_2,D_2)$ may be decomposed into a sequence of maps
\begin{align*}
(X_1,D_1) = (X^{(0)},D^{(0)})\dashrightarrow\cdots\dashrightarrow (X^{(N)},D^{(N)})=(X_2,D_2)
\end{align*}
where each of the arrows are blow-ups or blow-downs along smooth centers which have normal crossings with respect to the components of $D^{(j)}$. It therefore suffices to define $Ell(X,D;z,\tau)$ for smooth pairs $(X,D)$, where $D$ is a simple normal crossing divisor, and prove that $Ell(X,D;z,\tau)$ is functorial with respect to blow-ups along smooth centers which have normal crossings with respect to the components of $D$. This procedure has been carried out successfully by Borisov-Libgober in \cite{BL_Sing}, and by Chin-Lung Wang in \cite{CLW}. They define $Ell(X,D;z,\tau)$ by the formula:
\begin{align}\label{ell pair}
\int_X\prod_j\frac{x_j\vartheta(\twopi{x_j}-z,\tau)}{\vartheta(\twopi{x_j},\tau)}
\prod_i \frac{\vartheta(\twopi{D_i}-(a_i+1)z,\tau)\vartheta(z,\tau)}{\vartheta(\twopi{D_i}-z,\tau)\vartheta((a_i+1)z,\tau)}
\end{align} 
In the above expression, $x_j$ denote the formal Chern roots of $TX$ and $D_i$ denote the first Chern classes of the components $D_i$ of $D$ with coefficients $a_i(X,D)$. Note that since $\vartheta(0,\tau)=0$, the above expression only makes sense for $a_i\neq -1$. Naturally, this places some restrictions on the types of singularities allowed in the definition of $Ell(Y;z,\tau)$. For example, at the very least $Y$ must possess a log resolution $(X,D)\rightarrow (Y,0)$ such that none of the discrepancy coefficients $a_i(X,D)$ are equal to $-1$. In order to ensure that $Ell(Y;z,\tau)$ does not depend on our choice of a log resolution $(X,D)$, we actually must require that the discrepancy coefficients $a_i(X,D) > -1$. To see why, suppose that $(X_1,D_1)$ and $(X_2,D_2)$ are two log resolutions of $Y$ with discrepancy coefficients $a_i(X_j,D_j) \neq -1$. To prove that $Ell(X_1,D_1;z,\tau)=Ell(X_2,D_2;z,\tau)$, we must connect these two resolutions by a sequence of blow-ups and blow-downs, applying functoriality of the elliptic genus of pairs at each stage. But if some of the discrepancy coefficients $a_i(X_1,D_1)$ are greater than $-1$, and others less than $-1$, then after blowing up $X_1$, we may acquire discrepancy coefficients equal to $-1$. In this case, the elliptic genus of one of the intermediate pairs in the chain of blow-ups and blow-downs will be undefined, and consequently we will have no means of comparing the elliptic genera of $(X_1,D_1)$ and $(X_2,D_2)$. The only obvious way of avoiding this problem is to require $a_i(X_j,D_j) > -1$. This constraint is quite familiar to minimal model theorists; singular varieties $Y$ possessing this property are said to have \emph{log-terminal} singularities. 

Functoriality of the elliptic genus provides a nice explanation for the invariance of the elliptic genus under classical flops. For if $X_1$ and $X_2$ are related by a classical flop, then there exists a common resolution $f_i:X\rightarrow X_i$ with $f_1^*K_{X_1}= f_2^*K_{X_2}$. Two varieties related in this way are said to be \emph{$K$-equivalent}. One therefore discovers from this approach that the fundamental relation leaving the elliptic genus invariant is not flopping but $K$-equivalence.

Borisov-Libgober and Chin-Lung Wang's original proof of functoriality of the elliptic genus is based on an explicit calculation of the push-forward $f_*$ of the integrand in (\ref{ell pair}), where $f:(X,D)\rightarrow (X_0,D_0)$ is a blow-down. The obstruction to this push-forward giving the correct integrand on $X_0$ is given by an elliptic function with values in $H^*(X_0)$. One can then use basic elliptic function theory to show that this function vanishes. In what follows, we will sketch a different proof, similar to the one in \cite{EllipALE}, that makes use of the rigidity properties of the elliptic genus. This approach has several advantages: the first is that the proof can be easily generalized to more exotic versions of elliptic genera, such as the character-valued elliptic genus for orbifolds. Though the original proofs could be adapted to this situation, their implementation in the most general setting is cumbersome. Another advantage is that some variation of this approach appears to be useful for studying elliptic genera for varieties with non-log-terminal singularities. We will have more to say on this in the following section. Recall though that the rigidity of the elliptic genus for $SU$-manifolds was the key step in Totaro's proof of the invariance of elliptic genera under classical flops. It is therefore reasonable to expect rigidity phenomena to play a useful role in the study of elliptic genera of singular varieties.

Proceeding with the proof, we let $X$ be a smooth variety and $D=\sum a_i D_i$ a simple normal crossing divisor on $X$. Let $f:\Bl{X}\rightarrow X$ be the blow-up along a smooth subvariety which has normal crossings with respect to the components of $D$. We let $\Bl{D} =\sum a_i\Bl{D}_i + mE$ be the sum of the proper transforms of $D_i$ and the exceptional divisor $E$, whose coefficients are chosen so that $K_{\Bl{X}}-\Bl{D}=f^*(K_X-D)$.  

To avoid getting bogged down in technical details, let us assume that $f:\Bl{X}\rightarrow X$ is the blow-up at a single point $p = D_1\cap...\cap D_n$, and that $D_1,...,D_n$ are the only components of $D$. Then $T\Bl{X}$ is stably equivalent to $f^*TX\oplus\bigoplus_{i=1}^n \mathcal{O}(\Bl{D}_i)\oplus\mathcal{O}(E)$. The proof of the blow-up formula for the elliptic genus then amounts to proving that
\begin{align*}
&\int_{\Bl{X}}f^*\bigg\{\prod_{T'X}\ellip{x_j}\bigg\}\prod_{i=1}^{n}\frac{\twopi{\Bl{D}_i}\vartheta(\twopi{\Bl{D}_i}-(a_i+1)z)\vartheta'(0)}{\vartheta(\twopi{\Bl{D}_i})\vartheta(-(a_i+1)z)}\times\\
&\frac{\twopi{E}\vartheta(\twopi{E}-(m+1)z)\vartheta'(0)}{\vartheta(\twopi{E})\vartheta(-(m+1)z)} =\\
&\int_X\prod_{T'X}\ellip{x_j}\prod_{i=1}^n\frac{\twopi{{D}_i}\vartheta(\twopi{{D}_i}-(a_i+1)z)\vartheta'(0)}{\vartheta(\twopi{{D}_i})\vartheta(-(a_i+1)z)}
\end{align*}
Here, for ease of exposition, we have omitted the dependence of $\vartheta$ on $\tau$. Note that $\Bl{D}_i = f^*D_i-E$ in the above expression. Thus, if we expand both sides in the variables $f^*D_i, E$, and $D_i$, the blow-up fomula is easily seen to hold for integrals of Chern and divisor data not involving $E$. Note however that in a neighborhood of $E$, $(\Bl{X},\Bl{D})$ has the exact same structure as the blow-up of $\C^n$ at the origin, with the divisors $\Bl{D}_1,...,\Bl{D}_n$ corresponding to the proper transforms of the coordinate hyperplanes of $\C^n$. For the purpose of proving the blow-up formula, we may therefore assume that $X\cong (\Proj^1)^n$ and that $\Bl{X}$ is the blow-up of $X$ at $[0:1]\times\cdots\times[0:1]$. Viewed as a toric variety, $X$ is defined by the fan $\Sigma \subset N\otimes \R$ with $1$-dim rays $\R (\pm e_1),...,\R (\pm e_n)$, where $e_1,...,e_n$ form an integral basis for the lattice $N$. The fan $\Bl{\Sigma}$ of $\Bl{X}$ is obtained from $\Sigma$ by adding the ray $\R(e_1+...+e_n)$. The divisors $D_i\subset X$ correspond to the rays $\R e_i$ in $\Sigma$; and the divisors $\Bl{D}_i$ and $E$ correspond to the rays $\R e_i$ and $\R(e_1+...+e_n)$ in $\Bl{\Sigma}$. Using the fact that the tangent bundle of smooth toric variety with $T_\C$-invariant divisors $D_j$, $j=1,...,\ell$, is stably equivalent to $\bigoplus_{j=1}^\ell\mathcal{O}(D_j)$, the blow-up formula for $X$ reduces to proving:
\begin{align*}
&\int_{\Bl{X}}\prod_{k=1}^{n+1}\frac{\twopi{\Bl{D}_k}\vartheta(\twopi{\Bl{D}_k}-(a_k+1)z)\vartheta'(0)}{\vartheta(\twopi{\Bl{D}_k})\vartheta(-(a_k+1)z)}
\prod_{k=1}^{n}\frac{\twopi{\Bl{D}_{-k}}\vartheta(\twopi{\Bl{D}_{-k}}-(a_{-k}+1)z)\vartheta'(0)}{\vartheta(\twopi{\Bl{D}_{-k}})\vartheta(-(a_{-k}+1)z)}=\\
&\int_X \prod_{j=1}^n\frac{\twopi{{D}_j}\vartheta(\twopi{{D}_j}-(a_j+1)z)\vartheta'(0)}{\vartheta(\twopi{{D}_j})\vartheta(-(a_j+1)z)}
\prod_{j=1}^n\frac{\twopi{{D}_{-j}}\vartheta(\twopi{{D}_{-j}}-(a_{-j}+1)z)\vartheta'(0)}{\vartheta(\twopi{{D}_{-j}})\vartheta(-(a_{-j}+1)z)}
\end{align*}
In the above formula, $D_{-j}$ denote the $T_\C$-divisors on $X$ corresponding to the $1$-dim rays $\R(-e_j)$, with coefficients $a_{-j}=0$. $\Bl{D}_{-j}$ denote their proper transforms, which are simply given by $f^*D_{-j}$, since $D_{-j}$ are defined away from the blow-up locus. For ease of exposition, we also let $\Bl{D}_{n+1} = E$, with $a_{n+1}=m$. 

Now the crucial observation is that in the above formula, $\hbox{RHS}-\hbox{LHS}$ is independent of the coefficients $a_{-j}$. For since $\Bl{D}_{-j}$ are disjoint from $E$, any divisor intersection data involving $\Bl{D}_{-j}$ will be unchanged after formally setting $E=0$. Therefore, the parts of $\hbox{RHS}-\hbox{LHS}$ depending $a_{-j}$ will be unchanged after setting $E=0$. But formally letting $E=0$ clearly gives $\hbox{RHS}=\hbox{LHS}$. Consequently, $\hbox{RHS}-\hbox{LHS}$ depends only on $a_1,...,a_n$.

Let us therefore define $a_{-j}$ so that $(1+a_{-j}) = -(1+a_j)$. As discussed in the previous section, the set of coefficients $(1+a_{\pm j})$ assigned to the rays $\R(\pm e_j)$ give rise to a piece-wise linear functional $g = g_{1+a_i,1+a_{-i}}$ on the fan $\Sigma$. In fact, $g$ is simply the global linear functional which maps the basis vectors $e_i$ to $(1+a_i)$. As $g\in \hbox{Hom}(N,\Z)$, it also defines a global linear functional on $\Bl{\Sigma}$, taking the value $\sum_{i=1}^n (1+a_i)$ on $e_1+...+e_n$. Now by the discrepancy formula for blow-ups, $\sum_{i=1}^n (1+a_i) = (1+m)$. We see from this that the piece-wise linear functional on $\Bl{\Sigma}$ defined by assigning the coefficients $(1+a_{\pm j})$ to $\R(\pm e_j)$ and $(1+m)$ to $\R(e_1+...+e_n)$ corresponds to this same global linear functional $g$.

It follows that bundles $\mathcal{O}(e_{1}+...+e_{n})^{1+m}\otimes_{i=1}^n\mathcal{O}(e_i)^{1+a_i}\otimes\mathcal{O}(-e_i)^{1+a_{-i}}$ and $\otimes_{i=1}^n\mathcal{O}(e_i)^{1+a_i}\otimes\mathcal{O}(-e_i)^{1+a_{-i}}$ are trivial as $\Q$-line bundles on $\Bl{X}$ and ${X}$, respectively. Consequently, 
\begin{align*}
\hbox{Ind }\dbar\otimes\xi(1+a_i,1+m,1+a_{-i}) = \hbox{Ind }\dbar\otimes\xi(1+a_i,1+a_{-i}) = 0.
\end{align*}
But, up to a normalization factor, $\hbox{Ind }\dbar\otimes\xi(1+a_i,1+m,1+a_{-i}) = \hbox{RHS}$ and $\hbox{Ind }\dbar\otimes\xi(1+a_i,1+a_{-i}) = \hbox{LHS}$ for the given new values of $a_{-i}$. Thus, $\hbox{RHS}=\hbox{LHS}$ for $(1+a_{-i}) = -(1+a_i)$, and therefore also for $a_{-i} = 0$.

This completes the proof of the blow-up formula for the case where the blow-up locus is a single point. For completeness, let us outline the case for the blow-up along a subvariety $Z$ with normal crossings with respect to the components of $D$. This case is handled in much the same way, the only difference being that instead of reducing to the situation where $X$ is toric, we instead reduce to the case where $X$ is a toric fibration, fibered over the blow-up locus $Z$. Namely, by deformation to the normal cone, we may assume that $X = \Proj(M\oplus 1)\times\Proj(L_1\oplus 1)\times\cdots\times\Proj(L_r\oplus 1)$. Here, for the components $D_i$ intersecting $Z$, $L_i = \mathcal{O}(D_i)|_Z$ and $M$ is the quotient of $N_{Z/X}$ by $\oplus L_i$. The product $\times$ is the fiber product of the corresponding projective bundles over $Z$. We now view $D_i$ as the divisors given by the zero sections of the bundles $L_i$. Moreover, the zero sections of $L_i$ and $M$ together define a copy of $Z$ in $\Proj(M\oplus 1)\times\Proj(L_1\oplus 1)\times\cdots\times\Proj(L_r\oplus 1)$ with the same normal bundle $N_{Z/X}$ as in the original space. We let $\Bl{X}$ be the blow-up along this copy of $Z$. The proof of the blow-up formula then follows the same reasoning as in the toric case, where we now make use of the rigidity of fiber-wise analogues of the operators $\dbar\otimes\xi(\vec{a})$. For example, let us examine how to generalize the bundle $\xi(1+a_i,1+a_{-i})$ on $(\Proj^1)^n$ to the fibration $X$.

For each fibration $\pi_i: \Proj(L_i\oplus 1)\rightarrow Z$, we have the following exact sequence of tautological bundles
\begin{align*}
0\to S_i\rightarrow \pi^*(L_i\oplus 1)\rightarrow Q_i\to 0
\end{align*}
The vertical tangent bundle to $\Proj(L_i\oplus 1)$ is stably equivalent to the direct sum of hyperplane bundles $H_i\oplus H_{-i}$, where $H_i = \hbox{Hom}(\pi_i^*L_i,S_i)$ and $H_{-i}=\hbox{Hom}(1,S_i)$. Similarly, the vertical tangent bundle to the fibration $\pi:\Proj(M\oplus 1)\rightarrow Z$, with tautological bundle $S$ is stably equivalent to the direct sum $V\oplus H$ where $V = \hbox{Hom}(\pi^*M,S)$ and $H=\hbox{Hom}(1,S)$. All of these bundles extend naturally to the whole fibration $X$. Recall that if $\alpha_i = -\alpha_{-i}$, then $\dbar\otimes\xi(\alpha_i,\alpha_{-i})$ defines a elliptic operator on $(\Proj^1)^n$ with vanishing equivariant index (note that for convenience of notation we have defined $\alpha_i = 1+a_i$). For the fibration $X$, we replace $\dbar\otimes\xi(\alpha_i,\alpha_{-i})$ by the following fiber-wise analogue:
\begin{align*}
\dbar\otimes\bigotimes_{i=\pm 1}^{\pm r}&\bigotimes_{n=1}^{\infty}\Lambda_{-y^{\alpha_i}q^{n-1}}H_i^*\otimes\Lambda_{y^{-\alpha_i}q^n}H_i\otimes\bigotimes_{m=1}^{\infty}S_{q^m}H_i^*\otimes S_{q^m}H_i\otimes\\
&\bigotimes_{n=1}^{\infty}\Lambda_{-yq^{n-1}}V^*\otimes\Lambda_{y^{-1}q^n}V\otimes\bigotimes_{m=1}^{\infty}S_{q^m}V^*\otimes S_{q^m}V\otimes\\
&\bigotimes_{n=1}^{\infty}\Lambda_{-y^{-d-1}q^{n-1}}H^*\otimes\Lambda_{y^{d+1}q^n}H\otimes\bigotimes_{m=1}^{\infty}S_{q^m}H^*\otimes S_{q^m}H
\end{align*}
Here $d = \hbox{rank}(M)$. By performing a fiber integration over $X$, one can show that the rigidity of this operator with respect to the obvious torus action on the fibers follows directly from the rigidity results obtained for $\dbar\otimes\xi(\alpha_i,\alpha_{-i})$. Analogously, there exists a natural generalization of the operator $\dbar\otimes\xi(1+a_i,1+m,1+a_{-i})$ to a rigid operator on the fibration $\Bl{X}$. We therefore see that the blow-up formula for the elliptic genus is in all cases a consequence of rigidity phenomena on toric varieties.

Before moving on, we make a simple observation which will prove convenient in the next section. Let $X$ be a smooth toric variety with toric divisors $D_i$. Since $TX$ is stably equivalent to $\bigoplus_{i=1}^\ell\mathcal{O}(D_i)$, the elliptic genus of the pair $(X,\sum a_i D_i)$ is equal to the index of the operator $\dbar\otimes \xi(a_1+1,...,a_\ell+1)$, up to a normalization factor. Moreoever, the condtion that $\bigotimes_{i=1}^{\ell}\mathcal{O}(D_i)^{a_i+1}$ is trivial is equivalent to the condition that $K_X-\sum a_i D_i = 0$ as a Cartier divisor. In this case, we say that $(X,\sum a_i D_i)$ is a Calabi-Yau pair. Hence, a trivial consequence of the rigidity theorem for the elliptic genus of toric varieties is that $Ell(X,D;z,\tau)=0$ whenever $(X,D)$ is a toric Calabi-Yau pair.

\section{Beyond log-terminal singularities}
As observed above, Borisov-Libgober, and Chin-Lung Wang's approach to defining the elliptic genus of a singular variety $Y$ only appears to work when $Y$ has log-terminal singularities. This is due to the division by $\vartheta((a_i+1)z)$ in the formula for the elliptic genus of the pair $(X,D)$, where $(X,D)$ is a resolution of $Y$ with discrepancy coefficients $a_i(X,D)$. In pursuit of the broader question, ``for what class of singularities can we make sense of Chern data?", it is natural to ask whether log-terminality represents an essential constraint. In what follows, we will demonstrate that at the very least, the elliptic genus can be defined for all but a finite class of normal surface singularities.

Since the terms $\vartheta((a_i+1)z)$ do not involve any Chern data, the first thing one might try to do is simply throw these terms away in the definition of the elliptic genus of a pair. However, this approach is of little use since one would lose functoriality with respect to birational morphisms. As a second attempt, one could introduce a perturbation $a_i+\varepsilon b_i$ to each of the discrepancy coefficients $a_i$ of $D$, and take the limit as $\varepsilon\to 0$. Two obvious difficulties with this approach are $(1)$ the limit does not always exist, and $(2)$, even when the limit exists, it depends on the choice of the perturbation. Moreover, deciding on some fixed perturbation in advance (like letting all $b_i=1$) runs into problems if we hope to preserve functoriality.

To carry out this perturbation approach, we therefore require a distinguished class of peturbation divisors $\Delta(X,D) = \set{\sum\varepsilon b_i D_i}$ satisfying the following two properties: $(1)$ for every $D_\varepsilon\in \Delta(X,D)$, the limit as $\varepsilon\to 0$ of $Ell(X,D+D_\varepsilon;z,\tau)$ exists and is independent of the choice of $D_\varepsilon$; $(2)$ if $f:(\Bl{X},\Bl{D})\rightarrow (X,D)$ is a blow-up, then $f^*\Delta(X,D)\subset\Delta(\Bl{X},\Bl{D})$. Assuming we have found a set of perturbation divisors satisfying these properties, we could then define the elliptic genus of a singular variety $Y$ by the following procedure: Pick a log-resolution $(X,D)$ of $Y$, and choose $D_\varepsilon\in \Delta(X,D)$. Then define $Ell(Y;z,\tau) = \lim_{\varepsilon\to 0}Ell(X,D+D_\varepsilon;z,\tau)$. The important point is that if $f:(\Bl{X},\Bl{D})\rightarrow (X,D)$ is a blow-up, and $\Bl{D}_\varepsilon\in \Delta(\Bl{X},\Bl{D})$, then the answer we get for the elliptic genus of $Y$ is the same, regardless of whether we work with $(X,D+D_\varepsilon)$ or with $(\Bl{X},\Bl{D}+\Bl{D}_\varepsilon)$. To see why, note that $f^*(K_X-D-D_\varepsilon) = K_{\Bl{X}}-\Bl{D}-f^*D_\varepsilon$. Thus, by functoriality of the elliptic genus with respect to blow-ups, $Ell(X,D+D_\varepsilon;z,\tau)=Ell(\Bl{X},\Bl{D}+f^*D_\varepsilon;z,\tau)$. By property $(2)$, $f^*D_\varepsilon$ lies inside $\Delta(\Bl{X},\Bl{D})$. Hence, property $(1)$ of $\Delta(\Bl{X},\Bl{D})$ implies that $\lim_{\varepsilon\to 0}Ell(\Bl{X},\Bl{D}+\Bl{D}_\varepsilon;z,\tau)=\lim_{\varepsilon\to 0}Ell(\Bl{X},\Bl{D}+f^*D_\varepsilon;z,\tau)$.

For the case of complex surfaces, we have a natural candidate for $\Delta(X,D)$ satisfying the second property; namely the set
\begin{align*}
\set{\Delta_\varepsilon: \Delta_\varepsilon {D_i} = 0 \hbox{ for all }D_i\subset D \hbox{ with discrepancy coefficient} = -1}
\end{align*}
For if $\Delta_\varepsilon$ is in this set, and $\Bl{D}_i\subset \Bl{D}$ has coefficient equal to $-1$, then $f^*\Delta_\varepsilon\Bl{D}_i = \Delta_\varepsilon f_*\Bl{D}_i$. Now, either $\Bl{D}_i$ is the proper transform of a divisor with $-1$ discrepancy, or it is a component of the exceptional locus of $f$. In the former case, $\Delta_\varepsilon f_*\Bl{D}_i = 0$ by virtue $\Delta_\varepsilon$ belonging to the set $\Delta(X,D)$; in the latter case, $f_*\Bl{D}_i = 0$.

We still must verify that the $\varepsilon\to 0$ limit of $Ell(X,D+D_\varepsilon;z,\tau)$ is well-defined and independent of the choice of $D_\varepsilon\in \Delta(X,D)$ when $(X,D)$ is a resolution of a singular complex surface $Y$. Unfortunately, it is too much to ask that this property hold  for all normal surface singularities. Suppose, however, that $(X,D)$ is a log resolution of a normal surface $Y$ satisfying the following additional property: For every component $D_i\subset D$ with discrepancy coefficient $a_i(X,D) = -1$, $D_i \cong \Proj^1$ and $D_i$ intersects at most two other components of $D$ at a single point. In other words, we assume that the local geometry in a tubular neighborhood $U$ of $D_i$ is indistinguishable from a tubular neighborhood of a toric divisor. Note that since $D_i$ is an exceptional curve, the adjunction formula implies that $(X,D)|_U$ is a toric Calabi-Yau pair. Under this additional assumption, it turns out that $\lim_{\varepsilon\to 0}Ell(X,D+D_\varepsilon;z,\tau)$ exists and is independent of the choice of $D_\varepsilon\in \Delta(X,D)$. 

To see why the limit exists, note that $Ell(X,D+D_\varepsilon;z,\tau)$ is a meromorphic function of $\varepsilon$ with at most a simple pole at $\varepsilon=0$. Up to a normalization factor, the residue of $Ell(X,D+D_\varepsilon;z,\tau)$ at $\varepsilon=0$ corresponds to $\sum_{a_i(X,D)=-1}Ell(D_i,D+D_\varepsilon|_{D_i};z,\tau)$. By adjunction, $(D_i,D+D_\varepsilon|_{D_i})$ are all toric Calabi-Yau pairs, and consequently, the residue of $Ell(X,D+D_\varepsilon;z,\tau)$ vanishes by the rigidity theorems discussed in the previous section.

It remains to check that this limit is independent of the choice of $D_\varepsilon\in \Delta(X,D)$. Suppose then that $D_\varepsilon, D'_\varepsilon$ are two possible perturbation divisors. Since the $\varepsilon\to 0$ limit of $Ell(X,D+D_\varepsilon;z,\tau)-Ell(X,D+D'_\varepsilon;z,\tau)$ depends only on the local geometry near the divisor components $D_i$ with $a_i(X,D)=-1$, we may reduce the calculation to the case where $(X,D)$ is a toric variety. Moreover, since $(X,D)$ is Calabi-Yau in the tubular neighbhorhoods $U_i$ of the above divisor components, we may further reduce to the situation where $(X,D)$ is a Calabi-Yau pair. By definition, $D_\varepsilon$ and $D'_\varepsilon$ are trivial over $U_i$ and we may extend them to trivial divisors over all of $X$ without affecting the $\varepsilon\to 0$ limit of $Ell(X,D+D_\varepsilon;z,\tau)$ or $Ell(X,D+D'_\varepsilon;z,\tau)$ . We have thus reduced the calculation to the case where $(X,D+D_\varepsilon)$ and $(X,D+D'_\varepsilon)$ are both toric Calabi-Yau pairs. The rigidity theorem for the elliptic genus in this case then implies that $Ell(X,D+D_\varepsilon;z,\tau) = Ell(X,D+D'_\varepsilon;z,\tau) = 0$ for all $\varepsilon$, which clearly implies that their limits are the same as $\varepsilon\to 0$.

Of course, the above discussion is  moot unless one can find a reasonably large class of surface singularities whose resolutions satisfy the additional criterion of being locally toric in a neighborhood of the exceptional curves with $-1$ discrepancies. Fortunately, as observed by Willem Veys \cite{Veys}, nearly all normal surface singularities satisfy this property. The sole exceptions consist of the normal surfaces with strictly log-canonical singularities. These are surfaces whose resolutions $(X,D)$ satisfy $a_i(X,D)\geq -1$, with some $a_i(X,D)=-1$. A well-known example is the surface singularity obtained by collapsing an elliptic curve to a point. For a complete classification of these singularities, see \cite{Aster}. Based on this observation, Veys used a limiting procedure similar to the one given here to define Batyrev's string-theoretic Hodge numbers for normal surfaces without strictly log-canonical singularities. 

Note that, for dimensionality reasons, the elliptic genus of a smooth surface is a coarser invariant than the surface's collective Hodge numbers. Nevertheless, the approach discussed here affords several advantages. First, the technique of applying the rigidity properties of toric Calabi-Yau pairs is easy to adapt to more complicated invariants, such as the character-valued elliptic genus and the elliptic genus of singular orbifolds. These are finer invariants than the ordinary elliptic genus which are not characterized entirely by Hodge numbers. Second, this approach provides some clues about how to define elliptic genera for higher-dimensional varieties whose singularities are not log-terminal. For example, a possible generalization of the locally toric structure we required of the $-1$ discrepancy curves is to demand that all $-1$ discrepancy divisors be toric varieties fibered over some smooth base. The analogue of property $(2)$ for $\Delta(X,D)$ in this case is that  $c_1(D_\varepsilon)=0$ when restricted to each fiber of a $-1$ discrepancy divisor.

\section{Further directions}
Singular Chern numbers constructed out of elliptic genera have an interesting interpretation when the singular variety is the quotient of a smooth variety $X$ by a finite group $G$. In this situation, quantum field theory on orbifolds gives rise to a definition for the elliptic genus of $X/G$ constructed entirely out of the orbifold data of $(X,G)$. This orbifold version of the elliptic genus turns out to be closely related to the singular elliptic genus of $X/G$: for example, when the $G$-action has no ramification divisor, the orbifold elliptic genus of $(X,G)$ equals the singular elliptic genus of $X/G$. This is a specific example, proven by Borisov and Libgober \cite{BL_McKay}, of a much larger interaction between representation theory and topology known as the McKay correspondence.

Note that the log-terminality constraint comes for free in this case, since the germs of quotients $\C^n/G$, where $G$ is a finite subgroup of $GL(n,\C)$ are always log-terminal. Suppose however that $X$ itself is singular. By following a procedure similar to the one discussed above for the elliptic genus, one can continue to define a singular analogue of the orbifold elliptic genus of $(X,G)$. At this point it is natural to ask whether the McKay correspondence continues to hold when we allow $X$ to have singularities. When $X$ has log-terminal singularities, this follows directly out of Borisov and Libgober's proof of the McKay correspondence. For more general singularities the answer to this question is not known, although the McKay correspondence has been verified for the cases discussed in the previous section: namely, when $X$ is a normal surface without strictly log-canonical singularities. See, for example, \cite{RW_SingMcKay}.

As we have seen, many of the techniques for studying elliptic genera in birational geometry can be traced back to some rigidity property of the elliptic genus. It is therefore not surprising that most of these techniques (such as functoriality of the elliptic genus of a divisor pair) work equally well for the character-valued elliptic genus. From Totaro's work, we know that the elliptic genus completely determines the collection of Chern numbers invariant under classical flops. An obvious question then is whether the analogous statement holds for equivariant Chern numbers. From the functoriality property of the character-valued elliptic genus, one easily verifies that the equivariant Chern numbers encoded by the character-valued elliptic genus are indeed invariant under equivariant flops. The more difficult question is whether all flop-invariant equivariant Chern numbers factor through the character-valued elliptic genus. It appears that some knowledge of the image of the character-valued elliptic genus over an equivariant cobordism ring must play a role in answering this question.

\end{document}